\documentclass[12pt]{article}

\usepackage{amsfonts,amsmath,amsthm,latexsym,amssymb,fullpage,amscd}
\usepackage[mathscr]{eucal}
\usepackage{tikz}

 at 16 true pt
 at 13  true pt
 at 12 true pt

\newcommand{\R}{\mathbb{R}}
\newtheorem{theorem}{Theorem}[section]
\newtheorem{lemma}[theorem]{Lemma}
\newtheorem{corollary}{Corollary}[theorem]
\newtheorem{definition}{Definition}[section]

\title{Oscillatory integrals and weighted gradient flows}

\author{Michael Greenblatt}
\date{\today}

\newcommand\blfootnote[1]{%
  \begingroup
  \renewcommand\thefootnote{}\footnote{#1}%
  \addtocounter{footnote}{-1}%
  \endgroup
}
\begin{document}
\maketitle
\begin{abstract} 

We investigate estimating scalar oscillatory integrals by integrating by parts in directions based on $(x_1 \partial_{x_1} f(x) ,..., x_n \partial_{x_n}f(x))$, where $f(x)$ is the 
phase function. We prove a
theorem which provides estimates that are uniform with respect to linear perturbations of the phase and investigate some consequences. When the phase
function is quasi-homogeneous the theorem gives estimates for the associated surface measure Fourier transforms that are generally not too far off from being sharp.
In addition, the theorem provides a new proof, up to endpoints, that the well-known oscillatory integral estimates of Varchenko when the Newton polyhedron of the phase function 
is nondegenerate extend to corresponding bounds for surface measure Fourier transforms when the index is less than $\frac{1}{2}$. A sharp version of this was originally proven in [G2].
\vskip 0.1 in


\end{abstract}

\blfootnote{This work was supported by a grant from the Simons Foundation.}

\section{Background and theorem statements.} We consider oscillatory integrals of the form
\[I(\lambda) = \int  e^{i\lambda f(x_1,...,x_n)} \phi(x_1,...,x_n)\,dx_1\,...\,dx_n\tag{1.1}\]
Here  $f(x)$ is a real analytic function defined on a bounded neighborhood $U$ of the origin, $\phi(x)$ is a $C^1$ real-valued function supported in $U$,
and $\lambda$ is a real parameter. Often one seeks estimates of the form $|I(\lambda)| \leq h(|\lambda|)$, where
$h$ is an appropriately decreasing function, such as a function of the form $C(1 + |\lambda|)^{-s}$ for $C, s > 0$.  We will always assume that $\nabla f(0) = 0$ to ensure we are in a nontrivial
situation. By subtracting a constant from $f$, without loss of generality we may also assume that $f(0) = 0$.

A canonical example of where oscillatory integrals $(1.1)$ show up is in the analysis of Fourier transforms of surface measures. If $S$ is a surface in $\R^{n+1}$ that is given by the graph of  
a real analytic $f(x_1,...,x_n)$ on a bounded neighborhood $U$ of the origin and $\phi(x)$ is a real-valued $C^1$ function supported in $U$, then the Fourier transform of the Euclidean surface measure on $S$ localized through $\phi(x)$, which we denote by $\mu$, is given by
\[\hat{\mu}(\lambda_1,...,\lambda_{n+1}) = \int  e^{-i\lambda_{n+1} f(x_1,...,x_n) - i \lambda_1 x_1 - ... - i \lambda_n x_n} \phi(x_1,...,x_n)\,dx_1\,...\,dx_n\tag{1.2}\]
We always assume the surface has been translated and rotated so that like before, $f(0) = 0$ and $\nabla f(0) = 0$. This time the goal is to find estimates of the form  
$|\hat{\mu}(\lambda) | \leq h(|\lambda|)$ for appropriate decreasing $h$. 
Since Fourier transforms of surface measures appear in a range of subjects including  restriction
problems, maximal averages,  lattice point discrepancy, and more, estimates of this form can help lead to developments in those subjects. We  refer to Chapter 8 of [S] 
for an introduction to the connections between surface measure Fourier transforms and restriction problems, and we refer to [G5] for background material on 
the connections between maximal averages, lattice point discrepancy, and surface measure Fourier transforms. The paper [G5] also describes various results on surface measure Fourier transforms. 

Since the phase function in $(1.2)$ is a linear perturbation of that of $(1.1)$, one approach to proving such bounds on $|\hat{\mu}(\lambda) |$ is to provide bounds on $|I(\lambda)|$ that
are uniform under linear perturbations of the phase. In this paper we will describe one method of doing that. 

To motivate what we will be doing, observe that one approach to the analysis of oscillatory integrals $(1.1)-(1.2)$ is to appropriately divide the domain of integration into curves, do an 
appropriate integration by parts on each curve, and then integrate the result in the remaining $n-1$ dimensions. The hope would be that if the curves are properly chosen 
then one could obtain desirable functions $h(|\lambda|)$ bounding the overall oscillatory integral. A clue on how to 
select these curves is given by how oscillatory integral decay is often connected to bounds on sublevel set measures. Namely, in many situations, the supremum of the $\epsilon$ for which 
$(1.1)$  satisfies a bound $|I(\lambda)| \leq C(1 + |\lambda|)^{-\epsilon}$ for some $ C> 0$ is the same as the supremum of the $\epsilon$ for which there is a constant $C> 0$ such
 that the following holds for all $s > 0$.
\[ m(\{x \in U: |f(x)| < s\}) \leq C s^{\epsilon} \tag{1.3}\]
Here $m$ is Lebesgue measure. This suggests that one might get good results if the curves on which one performs the integrations by parts  are perpendicular to the boundaries of these sublevel sets.
In other words, one might choose these curves to be tangent to $\nabla f$. This idea is further bolstered by the fact that $f(x)$ increases or decreases fastest in the direction of
the gradient, so that the phase oscillates fastest in the directions of such curves.

However the above is not the whole story, since even if the phase is oscillating quickly in directions tangent to a curve, if the second derivative of the phase is also large these quick oscillations 
might not have the desired effect; the curve might head rapidly into region where the phase has a stationary point. Thus to add some flexibility to our activities, we consider not just 
curves whose tangents are in the direction of $\nabla f$, but also in directions $(a_1(x) \partial_{x_1} f(x) ,...,a_n(x) \partial_{x_n}f(x))$  where the weights $a_i(x)$ are real analytic or even
 quotients of real analytic functions. 

Using directions of this form  have an additional advantage. If one replaces $f(x)$ by $f(x) + b \cdot x$ for some $b \in \R^n$, then
$a_i (x) {\partial f \over \partial x_i} $ becomes $a_i(x)({\partial f \over \partial x_i} + b_i)$. Using resolution of singularities one can show that, 
generally speaking, for the type of $a_i(x)$ under discussion if one has sublevel set measure estimates of the form  $m(\{x \in U: \sum_{i = 1}^n |a_i (x) {\partial f \over \partial x_i}| < s\}) \leq C s^{\epsilon}$ for some $C > 0$ and
$0 < \epsilon < 1$, then if $U$ is sufficiently small the same will hold when each $a_i (x) {\partial f \over \partial x_i}$ is replaced by $a_i(x) ({\partial f \over \partial x_i} + b_i)$. 
This will allow us to state our theorems in terms of the optimal $\epsilon$ for which sublevel set bounds of the form $(1.3)$ hold for certain functions of the form 
$\sum_{i = 1}^n|\widetilde{a_i}(x) {\partial f \over \partial x_i}|$ (where the $\widetilde{a_i}(x)$ are slightly different from the $a_i(x)$.) Hence we will have a way of estimating $|I(\lambda)|$ that is uniform under linear perturbations, providing a way of bounding $|\hat{\mu}(\lambda)|$. 

To motivate possible choices of the weight functions $a_i(x)$, we consider the case where $f(x)$ is a monomial $a x_1^{\alpha_1}...\,x_n^{\alpha_n}$. When $\alpha_i > 0$, the 
 effect of taking an $x_i$ derivative on $f(x)$ is to multiply it by a constant times $x_i^{-1}$. This suggests that for "balance", one might choose $a_i(x) = x_i$, so that the directions in 
which one integrates by parts are of the form $(x_1 \partial_{x_1} f(x) ,..., x_n \partial_{x_n}f(x))$. Our main theorem, Theorem 1.1, will be based on using such weights. 
We will then see in section 2.1 
that this theorem provides good bounds for $|\hat{\mu}(\lambda)|$ when $f(x)$ is a quasi-homogeneous polynomial. Then in section 2.2 we will that Theorem 1.1  can be used to
show, up to endpoints, that when $f(x)$ has nondegenerate Newton polyhedron in the sense of [V], the optimal estimates for $|I(\lambda)|$
in [V] extend to analogous bounds for surface measure Fourier transforms. This was earlier shown (including endpoints) in [G2].

Other choices of $a_i(x)$ also lead to results that might be of some interest, but for simplicity of exposition we are only focusing on the weights $a_i(x) = x_i$ in this paper.

We now come to our theorem whose proof is based on using the above weighted gradient flow. While it gives especially desirable results in the above situations, the theorem holds generally.

\begin{theorem}
Suppose $\phi(x)$ is $C^1$ and $f(x)$ is real analytic with $f(0) = 0$ and $\nabla f(0) =  0$. Suppose $W$ is a bounded neighborhood of the origin such that if $\epsilon > 0$ is such 
that for some $C > 0$ and all $s > 0$ one has the sublevel set measure estimate
\[ m(\{x \in W: \frac{ \sum_{i = 1}^n |x_i {\partial f \over \partial x_i}|}{\prod_{i=1}^n |x_i|}  < s\}) \leq C s^{\epsilon} \tag{1.4}\]
Then if $U \subset W$ is a sufficiently small ball centered at the origin, the following hold. 
\begin{enumerate}

\item For all $\delta < \frac{\epsilon}{\epsilon + 1}$ there is a constant $A$ such that  $|I(\lambda)| \leq A(1 + |\lambda|)^{-\delta}$ whenever $\phi$ is supported in $U$. Here $A$ 
depends on $f$, $\phi$, and $\delta$.

\item For all 
$\delta < \min(\frac{\epsilon}{\epsilon + 1}, \frac{1}{2})$ there is a constant $B$ such that one has $|\hat{\mu}(\lambda)| \leq B(1 + |\lambda|)^{-\delta}$ whenever $\phi$ is supported in 
$U$. Here $B$ depends on $f$, $\phi$, and $\delta$.

\end{enumerate}

\end{theorem}

An examination of the proof of Theorem 1.1 reveals that the condition that $\phi$ is $C^1$ can actually  be weakened to requiring that there is a constant $M$ for which $|\phi(x)| \leq M$
 and $|\partial_{x_i} \phi(x)| \leq \frac{M}{|x_i|}$ for each $i$. On the other hand, if $\phi$ is smooth, 
it turns out that if $U$ is a small enough neighborhood of the origin, if $\epsilon$ is such that the sublevel set bounds $(1.3)$ hold for some $C$, then one automatically has
$|I(\lambda)| \leq C'(1 + |\lambda|)^{-\epsilon}$ for some $C'$ depending on $f$, $\phi$, and $U$. If in addition $\phi(x)$ is nonnegative with $\phi(0) > 0$, then the supremum of the $\epsilon$ for which such an estimate $|I(\lambda)| \leq C'(1 + |\lambda|)^{-\epsilon}$ holds is actually equal to the supremum of the $\epsilon$ for which $(1.3)$ holds, unless the former supremum is a negative integer. We refer to chapter 6 of [AGuV] for more information about these matters. This leads to the following corollary to Theorem 1.1, which may be of interest in its own right.

\begin{corollary}

If $W$ is a bounded neighborhood of the origin and $\epsilon_1$ denotes the supremum of the $\epsilon$ for which $(1.4)$ holds for some $C$, then if $U \subset W$ is a sufficiently small
neighborhood of the origin and $\epsilon_2$ denotes the supremum of the $\epsilon$ for which $(1.3)$ holds for some $C$, then  $\frac{\epsilon_1}{\epsilon_1 + 1} \leq \epsilon_2$. 

\end{corollary}

Although we won't prove it here, using resolution of singularities one can show that  $\epsilon_1$ is independent of $W$ and $\epsilon_2$ is independent of $U$ if
 $W$ and $U$ are sufficiently small neighborhoods of the origin, so that one can take $U = W$ in Corollary 1.1.1.

There has been quite a bit of work done on scalar oscillatory integrals of the form $(1.1)$. In addition to [V], some notable examples include the papers [BaGuZhZo] [CaCWr] [Gre] [Gr] [PhStS].
For the surface measure Fourier transforms $(1.2)$ there has also been a lot of work done, in part due to their connections with maximal averages, lattice point discrepancy, and other areas.
Much of the effort in this area has focused on either convex surfaces or the two-dimensional case. We mention the references [BNW] [BakMVaW] [BrHoI] [NaSeW] for the convex situation, and
[IkKeMu] [IkMu] for the two-dimensional situation. 

Gradient flows for scalar oscillatory integrals are often used in fields such as physics when putting the method of steepest descent into effect.
In addition, gradient flows appear in various areas of mathematics, including partial differential equations, optimization, and more applied fields like computer vision and machine learning.
The author does not know of any specific connection between this paper and the work in the above subjects, but it might be a direction worth exploring.

\section{Consequences of Theorem 1.1.}

\subsection {Quasi-homogeneous functions.}

A polynomial $p(x_1,...,x_n) = \sum_{\alpha} c_{\alpha} x^{\alpha}$ is said to be quasi-homogeneous if there are positive rational numbers $k_1,...,k_n$ such for any $\alpha =
(\alpha_1,...,\alpha_n)$ for which $c_{\alpha} \neq 0$, one has $\sum_{i=1}^n k_i \alpha_i = 1$. An equivalent statement is that $p(t^{k_1} x_1,...,t^{k_n} x_n) = t p(x_1,...,x_n)$ 
for all $(x_1,...,x_n) \in \R^n$ and all $t > 0$, and this definition extends the former definition to non-polynomials. A canonical example of a quasi-homogeneous polynomial is 
$x_1^{l_1} + ... + x_n^{l_n}$ where one has $k_i = {1 \over l_i}$ for each $i$.

Note that if $f(x)$ is a quasi-homogeneous polynomial, the function  $\sum_{i = 1}^n |x_i {\partial f \over \partial x_i}|$ appearing in Theorem 1.1
 is also quasi-homogeneous, with the same $(k_1,...,k_n)$.
Given the nature of the statement of Theorem 1.1 it makes sense that we would want to understand the growth rate of the measure of the sublevel sets of quasi-homogeneous functions.
 Suppose $f(x)$ is a quasi-homogeneous function and $k_1,...,k_n$ are as above. Let $V$ be a bounded neighborhood of the origin and let $\epsilon_0$ be defined by the supremum of the 
$\epsilon$ such that there is a constant $C > 0$ such that 
\[ m(\{x \in V: |f(x)| < s\}) \leq C s^{\epsilon} \tag{2.1}\]
Here as before $m$ denotes Lebesgue measure, and we will always work with functions for which $\epsilon_0 > 0$.  Note that an equivalent definition of $\epsilon_0$ is the supremum
 of the $\epsilon$ such that 
\[\int_V |f(x)|^{-\epsilon}\,dx < \infty \tag{2.2}\]
By the quasihomogeneity of $f$, the number $\epsilon_0$ is independent of $V$. In particular one may take $V$ to be the box $\{x: -1  < x_i < 1$ for all $i \}$, which we do henceforth.
If we change variables in $(2.2)$ from $x$ to $y$, where $x_1 = y_1$ and $x_i = ({\rm sgn\,} x_i) |y_i|^{k_i \over k_1}$ for $i > 1$, the integral in $(2.2)$ becomes a constant times
\[\int_V |\tilde{f}(y)|^{-\epsilon} \prod_{i = 1}^n |y_i|^{\frac{k_i}{k_1}- 1} \,dy \tag{2.3}\]
Here $\tilde{f}(y)$ is the function $f(x)$ in the $y$ coordinates, which has the key property that it is homogeneous of degree ${1 \over k_1}$. By conversion to polar coordinates, 
$(2.3)$ is finite when two requirements are met. First, we need that $-\epsilon{1 \over k_1} + \sum_{i = 1}^n (\frac{k_i}{k_1} - 1) > -n$. Secondly, we need that $|\tilde{f}(y)|^{-\epsilon}$
integrates to a finite value over the boundary sides of $V$, which is equivalent to $|f(x)|^{-\epsilon}$ integrating to a finite value over the boundary sides of $V$. The first condition 
translates into $\epsilon < \sum_{i=1}^n k_i$. Thus $\epsilon_0$ is the supremum of all $\epsilon$ such that $\epsilon < \sum_{i=1}^n k_i$ and such that for each $i$ and $a = 1,-1$ we have
\[ \int_{ \{ x:\, x_i = a, \, -1 < x_j < 1 {\rm\,\,for\,\,all\,\,} j \neq i\}} |f(x)|^{-\epsilon} \,dx_1,...,dx_{i-1}\,dx_{i+1},...,dx_n < \infty \tag{2.4}\]

We now investigate what Theorem 1.1 says in the quasi-homogeneous situation. Assume $f(x)$ is a quasi-homogeneous polynomial. Note that the function 
${\displaystyle \frac{ \sum_{i = 1}^n |x_i {\partial f \over \partial x_i}|}{\prod_{i=1}^n |x_i|}}$ appearing in Theorem 1.1 satisfies
\[ \frac{ \sum_{i = 1}^n |x_i {\partial f \over \partial x_i}|}{\prod_{i=1}^n |x_i|} \geq \sum_{i = 1}^n \bigg|x_i {\partial f \over \partial x_i}\bigg|\]
\[\geq C\bigg| \sum_{i = 1}^n k_ix_i {\partial f \over \partial x_i}\bigg|\]
\[ = C|f(x)| \tag{2.5}\]
Thus whenever we are in a situation where $(2.1)$ holds on some bounded neighborhood $V$ of the origin, the same will be true if we replace $|f(x)|$ by 
${\displaystyle \frac{ \sum_{i = 1}^n |x_i {\partial f \over \partial x_i}|}{\prod_{i=1}^n |x_i|}}$. In particular, this holds for every $\epsilon < \epsilon_0$. Thus part 2 of Theorem 1.1 says
that in the quasi-homogeneous case, one has $|\hat{\mu}(\lambda)| \leq B(1 + |\lambda|)^{-\delta}$ holds for all $\delta < \min({\epsilon_0 \over \epsilon_0 + 1}, {1 \over 2})$. So if 
$\epsilon_0 < 1$, whereas $\epsilon_0$ gives the supremal  exponent for the scalar oscillatory integral $I(\lambda)$ due to the connection between oscillatory integrals and sublevel set
measure growth, the supremal exponent for the surface measure Fourier transform is
 at least ${\epsilon_0 \over \epsilon_0 + 1}$. Hence the true exponent lies somewhere in the interval $[{\epsilon_0 \over \epsilon_0 + 1}, \epsilon_0]$. This can be a substantial improvement
over simply using stationary phase or the Van der Corput lemma along curves $(c_1 t^{k_1},...,c_n t^{k_n})$ as in [G4], where one can typically get an exponent no better than ${1 \over n + 1}$.

It might occur to one that since the estimates in $(2.5)$ are not that refined, we might be able to get better results by using more careful estimates than those of $(2.5)$. It turns out that
this often is the case if the condition $\epsilon < \sum_{i = 1}^n k_i$ is more stringent than the ones in $(2.4)$, so that $\epsilon_0 = \sum_{i = 1}^n k_i$. One can show that in 
many such situations, one will have that the exponent given by part 2) of  Theorem 1.1 is $\min(\sum_{i = 1}^n k_i, {1 \over 2})$, so that there is no reduction in the exponent
when it is less than ${1 \over 2}$. On the other hand, if the conditions of $(2.4)$ are more stringent than the  condition $\epsilon < \sum_{i = 1}^n k_i$, then often one can show that
Theorem 1.1 gives no better exponent than the $\min({\epsilon_0 \over \epsilon_0 + 1}, {1 \over 2})$ given above, while the true exponent can be as large as $\epsilon_0$.

\subsection{Estimates in terms of the Newton polyhedron.}

We first provide some relevant terminology.

\begin{definition} Let $f(x)$ be a smooth function defined on a neighborhood of the origin in 
$\R^n$, and
let $f(x) = \sum_{\alpha} f_{\alpha}x^{\alpha}$ denote the Taylor expansion of $f(x)$ at the origin.
For any $\alpha$ for which $f_{\alpha} \neq 0$, let $Q_{\alpha}$ be the octant $\{x \in \R^n: 
x_i \geq \alpha_i$ for all $i\}$. Then the {\it Newton polyhedron} $N(f)$ of $f(x)$ is defined to be 
the convex hull of all $Q_{\alpha}$. 
\end{definition} 

A Newton polyhedron may contain faces of dimensions zero through $n-1$ (a vertex is considered to be a compact face of dimension zero.)
These faces can be either compact or unbounded.  In this paper, as in earlier work like [G1] [G2] [V], an  important role is played by the following functions, defined for compact faces
 of the Newton polyhedron. 

\begin{definition}  Suppose $F$ is a compact face of $N(f)$. Then
if $f(x) = \sum_{\alpha} f_{\alpha}x^{\alpha}$ denotes the Taylor expansion of $f$ like above, we
define $f_F(x) = \sum_{\alpha \in F} f_{\alpha}x^{\alpha}$.
\end{definition} 

\begin{definition} The Newton polyhedron of $f(x)$ is said to be nondegenerate if for each compact face $F$ of $N(f)$, the function $\nabla f_F(x)$ is nonvanishing on $(\R - \{0\})^n$.
\end{definition} 

\begin{definition} Assume $N(f)$ is nonempty. Then the 
{\it Newton distance} $d(f)$ of $f(x)$ is defined to be $\inf \{t: (t,t,...,t,t) \in N(f)\}$.
\end{definition} 

\noindent The following is a well-known theorem of Varchenko [V]. We include the condition $|\lambda| \geq 2$ in the statement since the result is immediate for $|\lambda| < 2$ and we want
to avoid situations where $\ln |\lambda|$ is near zero. 

\begin{theorem} (Varchenko) Suppose $f(x)$ is real analytic on a neighborhood of the origin with $f(0) = 0$ such that $N(f)$ is nondegenerate. Let  $k$ denote
 the dimension of the face of $N(f)$ intersecting the line $x_1 = ... = x_n$ in its interior. There is a neighborhood $U$ of the origin such that if $\phi$ is supported in $U$ 
then the following hold .

\begin{enumerate}

\item There is a constant $C > 0$ depending on $f$ and $\phi$ such that $|I(\lambda)| \leq C|\lambda|^{-{1 \over d(f)}}(\ln|\lambda|)^
{n - 1 - k}$ for all $|\lambda| \geq 2$.

\item  If $\phi(x)$ is nonnegative, $\phi(0) 
> 0$, and $d(f) > 1$, then there is a $C' > 0$ depending on $f$ and $\phi$ such that $|I(\lambda)| \geq C'|\lambda|^{-{1 \over d(f)}}(\ln|\lambda|)^
{n - 1 - k}$ if $|\lambda|$ is sufficiently large.

\end{enumerate}
\end{theorem}

In [G2], among other things Theorem 2.1 was extended to surface measure Fourier transforms when $d(f)  > 2$. This follows from parts a) and  b) of the following consequence
 of Theorem 1.5 of [G2].

\begin{theorem} (Theorem 1.5 of [G2]).  Suppose $f(x)$ is real analytic on a neighborhood of the origin with $f(0) = 0$ and $\nabla f(0) =  0$. Let $k$ denote
 the dimension of the face of $N(f)$ intersecting the line $x_1 = ... = x_n$ in its interior. There
is a neighborhood $U$ of the origin such that if $\phi$ is supported in $U$ the following hold for $|\lambda|  \geq 2$, where $C$ denotes a constant depending on $f$ and $\phi$.

\noindent {\bf a)} If $d(f) < 2$, and each zero of each $f_F(x)$ on $(\R - \{0\})^n$ has order at most 2,
then there is a constant $C$ such that $|\hat{\mu}(\lambda)| \leq C |\lambda|^{-{1 \over 2}}$.

\noindent {\bf b)} If $d(f) \geq 2$ and each zero of each $f_F(x)$ on $(\R - \{0\})^n$ has order at most
$d(f)$, then there is a constant $C$ such that $|\hat{\mu}(\lambda)| \leq C|\lambda|^{-{1 \over d(f)}}(\ln |\lambda|)^
{n - k}$. If $d(f)$ is not an integer, the exponent $n - k$ can be improved to $n - 1 - k$. 

\noindent {\bf c)} If the maximum order $m$ of any zero of any $f_F(x)$ on $(\R - \{0\})^n$ 
satisfies $m > \max(d(f),2)$ then there is a constant $C$ such that $|\hat{\mu}(\lambda)| \leq C|\lambda|^{-{1 \over m}}$.

\end{theorem}

To be clear, the order of a zero of $f_F(x)$ at some $x_0$ here means the order of vanishing of the Taylor series of $f_F(x)$ at $x_0$. Also, the nondegeneracy condition in Theorem 2.1
 is equivalent to the statement that the zeroes of each $f_F(x)$ on $(\R - \{0\})^n$ are of order at most one. It may be worth pointing out that since the $f_F(x)$ here are quasi-homogeneous polynomials, their zeroes in $(\R - \{0\})^n$ are never isolated.

We now investigate what Theorem 1.1 says in the case where $f(x)$ has nondegenerate Newton polyhedron. The statement $(1.4)$ is equivalent to the statement that 
\[ m(\{x \in W: \frac{ \sum_{i = 1}^n (x_i {\partial f \over \partial x_i})^2}{\prod_{i=1}^n x_i^2}  < s\}) \leq C s^{\frac{\epsilon}{2}} \tag{2.6}\]
Let $g(x) = \sum_{i = 1}^n (x_i {\partial f \over \partial x_i})^2$. Then the Newton polyhedron $N(g)$ is the double $2N(f) = \{2x: x \in N(f)\}$, the faces $F'$ of $N(g)$ are the sets
$\{2x: x \in F\}$ for faces $F$ of $N(f)$. The statement that $N(f)$ is nondegenerate, namely that for each compact face $F$ of $N(f)$ the function $\nabla f_F(x)$ is 
nonvanishing on $(\R - \{0\})^n$, translates into the statement that each compact face $F'$ of $N(g)$, $g_{F'}(x)$  has no zeroes at all in $(\R - \{0\})^n$.

Suppose $d(f) > 1$, so that $d(g) > 2$. Let $h(x) = {\displaystyle \frac{ \sum_{i = 1}^n (x_i {\partial f \over \partial x_i})^2}{\prod_{i=1}^n x_i^2} = \frac{g(x)}{\prod_{i=1}^n x_i^2}}$.
Then one can define $N(h)$, $d(h)$, and $h_F(x)$ for faces $F$ of $N(h)$ analogously to Definitions 2.2-2.4. So we have $d(h) = d(g) - 2 = 2d(f) - 2 > 0$. Due to the analogous statement 
holding for $g(x)$, if $F$ is a compact face of $N(h)$ then the function $h_{F}(x)$ has no zeroes in $(\R - \{0\})^n$.

Note that the vertices of $N(h)$ 
may now have components as low as $-2$. Nonetheless many of the arguments of [G1] extend to $h(x)$, in particular the proof of Theorem 1.2 of [G1] which implies that since each
for each compact face $F$ of $N(h)$ the function $h_{F}(x)$ has no zeroes in $(\R - \{0\})^n$, equation
$(2.6)$ holds for all ${\epsilon \over 2} < {1 \over d(h)} = {1 \over 2d(f) - 2}$. This can also be  shown using toric resolution of singularities similarly to the arguments
in  [AGuV]. Consequently, $(1.4)$ holds for any  $\epsilon < {1 \over d(f) - 1}$, or equivalently when ${\displaystyle \frac{\epsilon}{\epsilon + 1} < \frac{\frac{1}{d(f) - 1}}{\frac{1}{d(f) - 1} + 1} =  {1 \over d(f)}}$.

 As a result, the second part of Theorem 1.1 says that $|\hat{\mu}(\lambda)| \leq B(1 + |\lambda|)^{-\delta}$
for each $\delta < \min({1 \over 2}, {1 \over d(f)})$. Up to endpoints, these are the estimates provided by Theorem 2.2. Such estimates are best possible since by taking 
$\lambda_k = 0$ for $k < n + 1$, one reduces to the oscillatory integral $I(\lambda)$ for which one has sharpness by [V] (part 2 of Theorem 2.1 here.) Thus we see that Theorem 1.1
provides another approach to proving such estimates, modulo endpoints.

\section{The proof of Theorem 1.1.}

\subsection {Some preliminary lemmas.}

We will make use of the following relatively easy lemma which follows from Lemma 3.2 of [G3].

\begin{lemma} (Lemma 3.2 of [G3]) Let $(E,\mu)$ be a finite measure space and suppose $g(x)$ is a measurable function on $E$ such that for some positive constants
 $C$ and $\delta$,
 for all $t  > 0$ one has $\mu(\{x \in E: |g(x)| < t\}) \leq Ct^{\delta}$. There is a constant $D_{\delta} > 0$
such that the following holds for all $M \neq 0$.
\begin{itemize}
\item If $\delta  < 1$, then $\int_E \min(1,  |M g|^{-1})\,d\mu < 
C D_{\delta}|M|^{-\delta}$
\item If $\delta = 1$, then $\int_E \min(1,|M g|^{-1})\,d\mu < 
C D_1(1 + \log_+ |M|)|M|^{-1} + \mu(E)|M|^{-1}$
\item If $\delta > 1$, then $\int_E \min(1,|M g|^{-1})\,d\mu < 
C(|M|^{-\delta} +  D_1 |M|^{-1})  + \mu(E)|M|^{-1}$
\end{itemize}
\end{lemma}

\noindent For the non-polynomial case, we will also need the following lemma from [G3].

\begin{lemma} (Corollary 2.1.2 of [G3]) Suppose $f_1(y_1,...,y_m)$,...,$f_l(y_1,...,y_m)$ are real analytic functions on a neighborhood of the origin, none identically zero. Then there
is an $m -1$ dimensional ball $B_{m-1}(0,\eta)$ and a positive integer $p$ such that for each $s_1,...,s_l$ and each $(y_1,...,y_{m-1}) \in B_{m-1}(0,\eta)$,
 the set $\{y_m: |y_m| < \eta $ and $f_i(y_1,...,y_m) < s_i$ for each $i\}$ consists of at most $p$ intervals.
\end{lemma}

\subsection {The beginning of the proof of Theorem 1.1.}

The size of the domain $U$ will be determined by our arguments; at certain junctures $U$ will have to be sufficiently small for the arguments to be valid. Also, we will always prove bounds 
of the form $C|\lambda|^{-\delta}$  for $|\lambda| \geq 2$ rather than $C(1 + |\lambda|)^{-\delta}$ for all $\lambda$ since the latter will always hold for $|\lambda| < 2$ simply 
by taking absolute values inside the integral and integrating. The exposition is somewhat easier if we prove estimates in the former form. In addition, we will always be bounding
 $|\hat{\mu}(\lambda)|$; bounds for $|I(\lambda)|$ will follow by setting $\lambda_k = 0$ for $k < n + 1$.

To start the proof, we observe that 
we may assume that $|(\lambda_1,...,\lambda_n)| < |\lambda_{n+1}|$, for if $U$ is sufficiently small, if $|(\lambda_1,...,\lambda_n)| \geq |\lambda_{n+1}|$ then the gradient of the
phase function is of magnitude at least $C|\lambda|$ and one may obtain a bound of $C'|\lambda|^{-1}$ simply by integrating by parts, better than the estimates we need.
Hence in our arguments
we will always assume that $|(\lambda_1,...,\lambda_n)| < |\lambda_{n+1}|$. This in particular implies that $|\lambda_{n+1}| \geq {1 \over \sqrt{2}} |\lambda|$. 

\noindent We next define the sets $U_i$ by
\[ U_i = \{x \in U: |x_i\partial_{x_i} f(x)| > |x_j\partial_{x_j} f(x)| {\rm\,\,for\,\,} j \neq i\} \tag{3.1}\]
As long as no two functions $x_j\partial_{x_j} f(x) $ are the same, up to a set of measure zero we will have $\cup_{i = 1}^n U_i = U$. In the rare event that two functions 
$x_j\partial_{x_j} f(x) $
 are in fact the same, we simply remove redundant $x_i \partial_{x_i} f(x) $ from the list when defining the $U_i$ in $(3.1)$, so that we always have $\cup_i U_i = U$ up to a set of
measure zero. The idea now is that since we are trying to integrate in directions based on the weighted gradient flow along $(x_1 \partial_{x_1} f(x) ,..., x_n \partial_{x_n}f(x))$, 
on each $U_i$ we will integrate by parts
in the $x_i$ direction, since in this direction $|x_i \partial_{x_i} f(x)|$ is at least $c|(x_1 \partial_{x_1} f(x) ,..., x_n \partial_{x_n}f(x))|$ for $c = n^{-{1 \over 2}}$.
 To this end, we correspondingly define the integrals $I_i(\lambda)$ by 
\[ I_i(\lambda) = \int_{U_i}  e^{-i\lambda_{n+1} f(x_1,...,x_n) - i \lambda_1 x_1 - ... - i \lambda_n x_n} \phi(x_1,...,x_n)\,dx_1\,...\,dx_n\tag{3.2}\]
Since $\cup_i U_i = U$ up to a set of measure zero, we have $\hat{\mu}(\lambda) = \sum_i I_i(\lambda)$, and in order to prove Theorem 1.1 it suffices to show that each $|I_i(\lambda)|$ is bounded
by the appropriate  of $A|\lambda|^{-\delta}$ or $B|\lambda|^{-\delta}$ as in the statement of the theorem.

 Let $P(x)$ be the function 
$ f(x_1,...,x_n) + {\lambda_1 \over \lambda_{n+1}} x_1 + ... +
{\lambda_n \over \lambda_{n+1}}x_n$, so that the phase function in $(3.2)$ is given by $\lambda_{n+1}P(x)$. The function $P(x)$ is sensible to focus on  here since $|\lambda_{n+1}| 
\sim |\lambda|$ and for our arguments it is helpful to view the phase function as a perturbation of the phase function when $\lambda_k = 0$ for all $k < n + 1$.

\noindent Where $\epsilon$ is such that $(1.4)$ holds, we write $U_i = D_1 \cup D_2$, where
\[D_1 = \{x \in U_i: |\partial_{x_i} P(x)| \leq |\lambda|^{-{1 \over \epsilon + 1}}\prod_{j \neq i} |x_j|\}\]
\[D_2 = \{x \in U_i: |\partial_{x_i} P(x)| >  |\lambda|^{-{1 \over \epsilon + 1}}\prod_{j \neq i} |x_j|\} \tag{3.3}\]
We correspondingly write $I_i(\lambda) = J_1(\lambda) + J_2(\lambda)$, where 
\[J_1(\lambda) = \int_{D_1}e^{-i\lambda_{n+1} f(x_1,...,x_n) - i \lambda_1 x_1 - ... - i \lambda_n x_n} \phi(x_1,...,x_n)\,dx_1\,...\,dx_n \tag{3.4a} \]
\[J_2(\lambda) = \int_{D_2}e^{-i\lambda_{n+1} f(x_1,...,x_n) - i \lambda_1 x_1 - ... - i \lambda_n x_n} \phi(x_1,...,x_n)\,dx_1\,...\,dx_n \tag{3.4b}\]
To bound $|J_1(\lambda)|$, we will simply take absolute values of the integrand and integrate in all variables, but it will take some effort to properly analyze the result.
 To bound $|J_2(\lambda)|$, we will perform integrations by parts in the $x_i$ variable.
 For this, we will need that each domain of integration in the $x_i$ variable consists of 
a number of intervals that is uniformly bounded. When $f(x)$ is a polynomial this is immediate, and for general real analytic $f(x)$ this follows from applying Lemma 3.2 as follows.

 We take the $y_m$ variable in that lemma to be the $x_i$ variable here, and the remaining $y_k$ variables to be the $x_j$ for $j \neq i$ along with two additional variables which
 we call $z_1$ and $z_2$. We take the $s_i$ of the lemma to be zero in all cases, and
we let the $f_i(y)$ of the lemma to be the functions $(x_j \partial_{x_j}f(x))^2 - (x_i\partial_{x_i}f(x))^2$ for $j \neq i$ along with the functions
$z_1\prod_{j \neq i} x_j^2 - (\partial_{x_i}f(x) + z_2)^2$  which we need for $z_1 = |\lambda|^{-{2 \over \epsilon + 1}}$
 and $z_2 = \frac{\lambda_i}{\lambda_{n+1}}$.
 
Although Lemma 3.2 is a local statement, by compactness one may extend the above application of Lemma 3.2 to the whole domain of integration here, so that there exists an $N$ such 
that for all $|\lambda| \geq 2$ and all $x_j$ for  $j \neq i$, the $x_i$ domain of integration of $J_2(\lambda)$ consists of at most $N$ intervals.  

We will bound $|J_1(\lambda)|$ and $|J_2(\lambda)|$ separately, starting with $|J_2(\lambda)|$.

\subsection  {The analysis of $J_2(\lambda)$.}

Since the exponents $\delta$ appearing in Theorem 1.1 satisfy $\delta < 1$, we can remove the regions where $|x_j| \leq |\lambda|^{-1}$ for some $j$ from the domains of integration when
bounding $|J_2(\lambda)|$. In other words, it suffices to bound $|J_2'(\lambda)|$, where $J_2'(\lambda)$ is given by 
\[ J_2'(\lambda) = \int_{\{x \in D_2:\, |x_j| > |\lambda|^{-1}{\rm\,\,for\,\,all\,\,} j\} }e^{-i\lambda_{n+1} f(x_1,...,x_n) - i \lambda_1 x_1 - ... - i \lambda_n x_n} \phi(x_1,...,x_n)\,dx_1\,...\,dx_n\tag{3.5}\]

For each $k = (k_1,...,k_n)$ such that the set $\{x: 2^{-k_j - 1} < |x_j| \leq 2^{-k_j}$ for all $j \}$ intersects the domain of integration of $(3.5)$, we define $J_{k,2}(\lambda)$ by 
\[ J_{k,2}(\lambda) = \int_{\{x \in D_2:\, |x_j| > |\lambda|^{-1}, \,2^{-k_j - 1} < |x_j| \,\leq\,2^{-k_j} {\rm\,\,for\,\,all\,\,}j \} }e^{-i\lambda_{n+1} f(x_1,...,x_n) - i \lambda_1 x_1 - ... - i \lambda_n x_n} \phi(x)\,dx_1\,...\,dx_n\tag{3.6}\]
Hence $J_2'(\lambda) = \sum_k J_{k,2}(\lambda) $.
Without loss of generality, we may assume $U$ is a subset of the unit ball, so that each $k_j \geq 0$ and there at most $C(\ln |\lambda|)^n$ terms $J_{k,2}(\lambda)$ to consider. 

In accordance with the weighted gradient flow concept, since we are on $U_i$, where $|x_i \partial_{x_i} f(x)| \geq n^{-{1 \over 2}}|(x_1 \partial_{x_1} f(x) ,..., x_n \partial_{x_n}f(x))|$, 
a given term $J_{k,2}(\lambda)$ will be analyzed using an integration by parts in the $x_i$ direction.

In $(3.6)$, we now integrate by parts in the $x_i$ variable over  in any of the at most $N$ intervals of integration for fixed other variables. Suppose $L = [l_1,l_2]$ is one such interval. 
Then the exponential appearing in $(3.6)$ can be written as $e^{-i\lambda_{n+1} P(x)}$, which we may write as $-i\lambda_{n+1}\partial_{x_i} 
P(x) (e^{-i\lambda_{n+1} P(x)}/ -i\lambda_{n+1}\partial_{x_i} P(x))$. Note that the denominator $-i\lambda_{n+1}\partial_{x_i} P(x)$ is never zero by the definition of $D_2$. We integrate
by parts in $x_i$ over $L$, integrating $-i\lambda_{n+1}\partial_{x_i} P(x) e^{-i\lambda_{n+1} P(x)}$ to $e^{-i\lambda_{n+1} P(x)}$ and differentiating the rest. Taking absolute values
in the result, we get the following, where $P_0(x_i)$ denotes $P(x)$ as a function of $x_i$ with all other variables fixed. 
\[ \bigg|\int_L  e^{-i\lambda_{n+1} f(x_1,...,x_n) - i \lambda_1 x_1 - ... - i \lambda_n x_n} \phi(x_1,...,x_n)\,dx_i\bigg|   \]
\[ \leq  C{1 \over |\lambda_{n+1}|} \bigg({1 \over |\partial_{x_i}P_0(l_1)|} + {1 \over |\partial_{x_i} P_0(l_2)|}+   
\int_L \frac{|\partial_{x_ix_i}^2 P_0(x_i)|}{(\partial_{x_i}P_0(x_i))^2}\,dx_i + \int_L \frac{1}{|\partial_{x_i}P_0(x_i)|}\,dx_i\bigg)  \tag{3.7} \]
We would like to integrate the first  integral in $(3.7)$ back to get terms similar to the endpoint terms of $(3.7)$. For this to work, we need that $\partial_{x_ix_i}^2 P_0(x_i)$ 
changes sign  at boundedly
many points in a given interval $L$. Since $\partial_{x_ix_i}^2 P_0(x_i)  = \partial_{x_ix_i}^2 f(x_1,...,x_n)$, we may show this by once again invoking Lemma 3.2, this time for one 
function $- (\partial_{x_ix_i}^2 f(x_1,...,x_n))^2$, letting the $y_m$ variable be $x_i$ and $y_1,...,y_{m-1}$ the remaining $x_j$ variables, and $s_1 = 0$. Again, due to compactness
it suffices to have the local result of Lemma 3.2 (or alternatively, one could just assume $U$ is a sufficiently small neighborhood of the origin.) In any event, there is an $N_0$ such that
$\partial_{x_ix_i}^2 P_0(x_i)$ changes sign at most $N_0$ times on any interval $L$. Consequently, one may write the integral in $(3.7)$ as the union of at most $N_0$ intervals on which
$\partial_{x_ix_i}^2 P_0(x_i)$ is either nonnegative or nonpositive, and integrate back to terms similar to the endpoint terms. Specifically, one may bound $(3.7)$ by
\[C{1 \over |\lambda_{n+1}|}\bigg({1 \over \partial_{x_i}P_0(l_1)|} + {1 \over |\partial_{x_i} P_0(l_2)|}+ \sum_l {1 \over |\partial_{x_i} P_0(j_l)|}
+ \int_L \frac{1}{|\partial_{x_i}P_0(x_i)|}\,dx_i\bigg)\tag{3.8} \]
Here there are at most $N_0$ terms in the sum of $(3.8)$. Since the domain of integration of $J_{k,2}$ consists of points where $|\partial_{x_i} P(x)| >  
|\lambda|^{-{1 \over \epsilon + 1}}\prod_{j \neq i} |x_j|$ and since $|\lambda_{n+1}| \geq {1 \over \sqrt{2}}|\lambda|$ in the situation at hand, we have that $(3.8)$ is bounded by
\[ C'|\lambda|^{-{\epsilon \over \epsilon + 1}}(\prod_{j \neq i} |x_j|)^{-1} \tag{3.9}\]
Equation $(3.9)$ bounds the $x_i$ integral over a single interval $L$. There are at most $N$ such intervals, so the overall integral of $(3.6)$ in the $x_i$ direction is also bounded
by a constant times $|\lambda|^{-{\epsilon \over \epsilon + 1}}(\prod_{j \neq i} |x_j|)^{-1}$. If one integrates in the remaining $n-1$ variables, keeping in mind that the size of the 
domain of integration in each $x_j$ direction for $j \neq i$ is bounded by a constant times $|x_j|$, we obtain that
\[|J_{k,2}(\lambda)| \leq C''|\lambda|^{-{\epsilon \over \epsilon + 1}}\tag{3.10}\]
Because there are at most a constant times $(\ln|\lambda|)^n$ possible values of $k$, the sum of all $|J_{k,2}(\lambda)|$ is at most 
$C_3(\ln|\lambda|)^n|\lambda|^{-{\epsilon \over \epsilon + 1}}$. This is less than $C_4|\lambda|^{-\delta}$ for any $\delta < \frac{\epsilon}{\epsilon + 1}$, so $|J_2(\lambda)| \leq 
\sum_k |J_{k,2}(\lambda)|$ will always satisfy the bounds needed for Theorem 1.1. Thus we may 
focus our attention henceforth on bounding $|J_1(\lambda)|$, which we must show is bounded by the appropriate quantity in Theorem 1.1.

\subsection {The analysis of $J_1(\lambda)$.}
 
Taking absolute values of the integrand in $(3.4a)$ and integrating (recalling the definition $(3.3)$ of $D_1$) leads to
\[ |J_1(\lambda)| \leq Cm(\{x \in U_i:  |x_i \partial_{x_i} P(x)| \leq |\lambda|^{-{1 \over \epsilon + 1}}\prod_{j=1}^n  |x_j|\})\tag{3.11}\]
Here as before $m$ denotes Lebesgue measure.
Suppose we are in the special case where $\lambda_k = 0$ for $k < n + 1$, the setting of part 1 of Theorem 1.1. Here $P(x) = f(x)$.
Since $U_i$ is the set of points where $|x_i \partial_{x_i}f(x)| > |x_j \partial_{x_j}f(x)|$ for all $j \neq i$, $(3.11)$ leads to
\[ |J_1(\lambda)| \leq Cm(\{x \in U: \sum_{j = 1}^n |x_j \partial_{x_j} f(x)| \leq n |\lambda|^{-{1 \over \epsilon + 1}} \prod_{j=1}^n  |x_j|\})\tag{3.12}\]
By the assumption $(1.4)$, we have that $(3.12)$ is bounded by a constant times $|\lambda|^{-{\epsilon \over \epsilon + 1}}$. This is better than what is needed for part 1 of Theorem 1.1, 
so we have now shown part 1 of Theorem 1.1.

We focus our attention on proving part 2 of the theorem.
We will show that if the neighborhood $U$ of the origin is sufficiently small, the right-hand side of $(3.11)$ is bounded by $C|\lambda|^{-\delta}$ whenever $\delta < \min(\frac{\epsilon}{\epsilon + 1}, \frac{1}{2})$. 

Since $\partial_{x_i} P(x) = \partial_{x_i} f(x) - c$ where $c = \frac{\lambda_i}{\lambda_{n+1}}$ satisfies $|c| \leq 1$, we can rewrite the measure on the
 right-hand side of $(3.11)$ as 
\[m(\{x \in U_i:  \bigg|{\partial_{x_i}f (x) - c \over \prod_{j \neq i}  x_j}\bigg| \leq |\lambda|^{-{1 \over \epsilon + 1}}\})\tag{3.13}\]
The idea behind bounding $(3.13)$ is that using resolution of singularities, one can show that if $\epsilon < 1$ then the measure in $(3.13)$ is bounded by the measure when $c = 0$, while
if $\epsilon \geq  1$ one at least has a uniform bound of $C|\lambda|^{-{1 \over \epsilon + 1}}$. In the former case, one uses the $c = 0$ case exactly above to get an overall bound of 
$C|\lambda|^{-{\epsilon \over \epsilon + 1}}$, while in the latter case since $\epsilon \geq 1$, $(1.4)$ holds if one replaces $\epsilon = 1$ and thus one can use the $\epsilon = 1$ case
to say that we have an overall bound of $C|\lambda|^{-{1 \over 2}}$. 

We use resolution of singularities, first proved by Hironaka [H1][H2], as follows. Let $W$ be as in Theorem 1.1. Assuming $U$ is a sufficiently small neighborhood of the origin, by resolution of 
singularities there exists an open $U' \subset W$ containing the closure $\overline{U}$ and a finite collection $\{g_l\}_{l=1}^L$ of  proper real analytic mappings $g_l: V_l \rightarrow U'$, 
where each $V_l$ is open and contains $\overline{g_l^{-1}(U)}$, with the following properties.

Each $g_l$ is one to one outside a set of measure zero.  Let $\{h_m(x)\}_{m=1}^M$ denote the list of functions
 consisting of each $\partial_{x_j} f(x)$, each nonzero difference $(x_j \partial_{x_j} f(x))^2  - (x_k \partial_{x_k} f(x))^2$, and each coordinate function $x_j$. Then on $V_l$
 each function $h_m \circ g_l(x)$ and the Jacobian determinant ${Jac\,}_l(x)$ of each $g_l(x)$  is of 
the form $a_{lm}(x) p_{lm}(x)$, where $p_{lm}(x)$ is a monomial and $a_{lm}(x)$ never vanishes; in fact one will have $|a_{lm}(x)| > \delta$ for some $\delta > 0$.
 For any bounded measurable function $F(x)$ on $U$ one has
\[ \bigg|\int_U F(x) \, dx \bigg| \leq \sum_{l=1}^L  \int_{g_l^{-1}(U)} |F \circ g_l (x) {Jac\,}_l(x)| \,dx \tag{3.14}\]
Because of our assumption that $\nabla f(0) = 0$, the monomials $p_{lm}(x)$ are never just the constant monomial $1$.

Let $m_0$ be such that
$h_{m_0}(x) = \partial_{x_i} f(x)$ for the $i$ being used above and let $j_0$ be any index such that such that $p_{lm_0}(x)$ contains $x_{j_0}$ to at least the first power. Then for some 
$d \neq 0$ one has
\[ \partial_{x_{j_0}} (a_{lm_0}(x)p_{lm_0}(x)) = \big(\partial_{x_{j_0}}a_{lm_0}(x) + {d \over x_{j_0}}a_{lm_0}(x)\big)p_{lm_0}(x) \tag{3.15}\]
As long as $U$ is a sufficiently small neighborhood of the origin, then the term ${d \over x_{j_0}}a_{lm_0}(x)$ 
will dominate $\partial_{x_{j_0}}a_{lm_0}(x)$ in absolute value for at least one choice of $j_0$ since we will be near enough to $g_l^{-1}(0)$. 
So we may work under the assumption that at each $x$ there is always some $j_0$ and some constant $d' > 0$ so that we have
\[ |\partial_{x_{j_0}} (a_{lm_0}(x)p_{lm_0}(x))| \geq {d' \over |x_{j_0}|}|a_{lm_0}(x)p_{lm_0}(x)| \tag{3.16}\]
In order to bound $(3.13)$, we will apply $(3.14)$, letting $F(x)$ be the characteristic function of $\{x \in U_i: {\displaystyle \bigg|{\partial_{x_i}f (x) - c \over \prod_{j \neq i}  x_j}\bigg|}\leq 
|\lambda|^{-{1 \over \epsilon + 1}}\}$, or equivalently the characteristic function of 
\[\{x \in U_i: c - |\lambda|^{-{1 \over \epsilon + 1}}\prod_{j \neq i}  |x_j| < \partial_{x_i}f (x) 
< c + |\lambda|^{-{1 \over \epsilon + 1}}\prod_{j \neq i}  |x_j|\} \tag{3.17}\]
Let $q_i(x) = \prod_{j \neq i}  x_j$. Then $q_i \circ g_l(x)$ is of the form $b_{il}(x)r_{il}(x)$ where $r_{il}(x)$ is a monomial and $|b_{il}(x)| > \delta$ for some positive $\delta$. Then 
$F \circ g_l(x)$ the characteristic function of 
\[\{x \in g_l^{-1}(U_i): c - |\lambda|^{-{1 \over \epsilon + 1}}|b_{il}(x)r_{il}(x)| < a_{lm_0}(x)p_{lm_0}(x)
< c+|\lambda|^{-{1 \over \epsilon + 1}}|b_{il}(x)r_{il}(x)|\} \tag{3.18}\] 
Denote the set in $(3.18)$ by $E_{il\lambda c}$. Then for the above $F(x)$, equation $(3.14)$ leads to
\[m(\{x \in U_i:  \bigg|{\partial_{x_i}f (x) - c \over \prod_{j \neq i}  x_j}\bigg| \leq |\lambda|^{-{1 \over \epsilon + 1}}\}) \leq \sum_{l=1}^L \int_{E_{il\lambda c}} |{Jac\,}_l(x)|\, dx \tag{3.19}\]

Thus we turn our attention to bounding a given term  $\int_{E_{il\lambda c}} |{Jac\,}_l(x)|\, dx$. We examine the intersection of the domain of integration $E_{il\lambda c}$ with a given dyadic rectangle which we denote by $R$. 
Because ${Jac\,}_l(x)$, $b_{il}(x)r_{il}(x)$, and $a_{lm_0}(x)p_{lm_0}(x)$ are comparable to monomials, on $E_{il\lambda c} \cap R$ the functions $|{Jac\,}_l(x)|$, $|b_{il}(x)r_{il}(x)|$, and 
$|a_{lm_0}(x)p_{lm_0}(x)|$ are within a constant factor of the functions
 $|{Jac\,}_l(x^*)|$, $|b_{il}(x^*)r_{il}(x^*)|$, and $|a_{lm_0}(x^*)p_{lm_0}(x^*)|$ respectively where $x^* \in E_{il\lambda c} \cap R$ is fixed.

As long as the neighborhood of $U$ is sufficiently small, not only does $(3.16)$ hold, but also for each $R$ there will necessarily be a single $j_0$ for 
which $(3.16)$ holds  throughout $R$, since we will be close enough to $g_l^{-1}(0)$ for this to be true. By $(3.16)$ and $(3.18)$,
the width of the $x_{j_0}$ cross section of $E_{il\lambda c} \cap R$ is bounded by 
\[C  |\lambda|^{-{1 \over \epsilon + 1}}|b_{il}(x^*)r_{il}(x^*) | \bigg|\frac{a_{lm_0}(x^*)p_{lm_0}(x^*)}{x_{j_0}^*}\bigg|^{-1} \tag{3.20}\]
The width of this cross section is also trivially bounded by $ x^*_{j_0}$. These can be combined by saying the width of this cross section is bounded by the following, where $R_0$ denotes 
the cross section of $R$ in the $x_{j_0}$ direction.
\[C \int_{R_0} \min(1,  |\lambda|^{-{1 \over \epsilon + 1}}|b_{il}(x^*)r_{il}(x^*)||a_{lm_0}(x^*)p_{lm_0}(x^*)|^{-1})\,dx_{j_0}\tag{3.21}\]
Integrating this in the remaining $n-1$ variables and inserting the Jacobian factor, we see that 
\[ \int_{E_{il\lambda c} \cap R} |{Jac\,}_l(x)| \, dx \leq C \int_R|{Jac\,}_l(x^*)|\min(1,  |\lambda|^{-{1 \over \epsilon + 1}}|b_{il}(x^*)r_{il}(x^*)||a_{lm_0}(x^*)p_{lm_0}(x^*)|^{-1})\,dx 
\tag{3.22}\]
Since the various factors in $(3.22)$ stay within a bounded factor on $R$, we may replace $(3.22)$ by
\[ \int_{E_{il\lambda c} \cap R} |{Jac\,}_l(x)| \, dx \leq C \int_R|{Jac\,}_l(x)|\min(1,  |\lambda|^{-{1 \over \epsilon + 1}}|b_{il}(x)r_{il}(x)||a_{lm_0}(x)p_{lm_0}(x)|^{-1})\,dx 
\tag{3.23}\]
We next add $(3.23)$ over all $R$ intersecting $g_l^{-1}(U_i)$. Letting $W_{il}$ denote the union of all $R$ intersecting $g_l^{-1}(U_i)$, we get
\[ \int_{E_{il\lambda c}} |{Jac\,}_l(x)| \, dx \leq C \int_{W_{il}}|{Jac\,}_l(x)|\min(1,  |\lambda|^{-{1 \over \epsilon + 1}}|b_{il}(x)r_{il}(x)||a_{lm_0}(x)p_{lm_0}(x)|^{-1})\,dx 
\tag{3.24}\]
By shrinking $U$ if necessary, we can always assume that $W_{il} \subset g_l^{-1}(U')$, where $U'$ is the open set containing $\overline{U}$ of the resolution of singularities procedure.
Thus $(3.24)$ becomes 
\[ \int_{E_{il\lambda c}} |{Jac\,}_l(x)| \, dx \leq C \int_{g_l^{-1}(U')}|{Jac\,}_l(x)|\min(1,  |\lambda|^{-{1 \over \epsilon + 1}}|b_{il}(x)r_{il}(x)||a_{lm_0}(x)p_{lm_0}(x)|^{-1})\,dx 
\tag{3.25}\]
Now going back into the original coordinates using the coordinate change map $g_l$, $(3.25)$ leads to 
\[ \int_{E_{il\lambda c}} |{Jac\,}_l(x)| \, dx \leq C \int_{U'}\min(1,  |\lambda|^{-{1 \over \epsilon + 1}}(\prod_{j \neq i}| x_j|)|\partial_{x_i} f(x) |^{-1})\,dx \tag{3.26}\]
We can refine $(3.26)$ as follows. Since each nonzero  difference $(x_j \partial_{x_j} f(x))^2  - (x_k \partial_{x_k} f(x))^2$ was monomialized, in particular $(x_i \partial_{x_i} f(x))^2 
 - (x_j \partial_{x_j} f(x))^2$ is monomialzed for each $j$. Thus in the blown up coordinates of $g_l^{-1}(U')$, the set of points where $|x_i \partial_{x_i} f(x)| > |x_j \partial_{x_j} f(x)| $ 
for a given $j$, which is the same as the set where
 $(x_i \partial_{x_i} f(x))^2  - (x_j \partial_{x_j} f(x))^2 > 0$, consists of the points in  $g_l^{-1}(U')$ where a certain monomial is positive. This consists of the set of points in 
$g_l^{-1}(U')$ that are in
a certain collection of octants. Thus in the blown up coordinates, the  set of points in  $g_l^{-1}(U')$  where $|x_i \partial_{x_i} f(x)| > |x_j \partial_{x_j} f(x)| $ for all $j$ are
 also the points in $g_l^{-1}(U')$ that are in certain octants. Furthermore, the above dyadic rectangles $R$ will be contained in these octants. Hence the relation $|x_i \partial_{x_i} f(x)|
 > |x_j \partial_{x_j} f(x)| $ for $j \neq i$ will still hold when we return the union of all such $R$ to the original coordinates as in $(3.26)$. Thus we may amend $(3.26)$ to 
\[ \int_{E_{il\lambda c}} |{Jac\,}_l(x)| \, dx \leq C \int_{\{x \in U':|x_i \partial_{x_i} f(x)| \, >\, |x_j \partial_{x_j} f(x)|{\rm\,\, for\,\,} j \neq i\}}\min(1,  |\lambda|^{-{1 \over \epsilon + 1}}(\prod_{j \neq i}| x_j|)|\partial_{x_i} f(x) |^{-1})\,dx \tag{3.27}\]
Since $(\prod_{j \neq i}| x_j|)|\partial_{x_i} f(x) |^{-1} = {\displaystyle \bigg({|x_i \partial_{x_i} f(x)| \over \prod_{j = 1}^n|x_j|}\bigg)^{-1}}$ and $|x_i \partial_{x_i} f(x)| \, >\, |x_j \partial_{x_j} f(x)|$ for all
 $j \neq i$ in the domain of integration of $(3.27)$, we may let $r(x) = {\displaystyle \frac{ \sum_{i = 1}^n |x_i {\partial f \over \partial x_i}|}{\prod_{i=1}^n |x_i|}}$ (the function in Theorem 1.1) 
and then $(3.27)$ implies
\[ \int_{E_{il\lambda c}} |{Jac\,}_l(x)| \, dx \leq C \int_{U'}\min(1,  |\lambda|^{-{1 \over \epsilon + 1}}|r(x)|^{-1})\,dx \tag{3.28}\]
Combining with $(3.19)$ leads to 
\[m(\{x \in U_i:  \bigg|{\partial_{x_i}f (x) - c \over \prod_{j \neq i}  x_j}\bigg| \leq |\lambda|^{-{1 \over \epsilon + 1}}\}) \leq C'\int_{U'}\min(1,  |\lambda|^{-{1 \over \epsilon + 1}}|r(x)|^{-1})\,dx \tag{3.29}\]
Since the left hand side of $(3.29)$ is the right-hand side of $(3.13)$, which by $(3.11)$ is an upper bound for $|J_1(\lambda)|$, we conclude that 
\[|J_1(\lambda)|  \leq C'\int_{U'}\min(1,  |\lambda|^{-{1 \over \epsilon + 1}}|r(x)|^{-1})\,dx \tag{3.30}\]
We now apply Lemma 3.1, using the definition $(1.4)$ of $\epsilon$ and the fact that $U' \subset W$. If $\epsilon < 1$, the first part of Lemma 3.1 gives that the right hand side of 
$(3.30)$ is bounded 
by a constant times $|\lambda|^{-{\epsilon \over \epsilon + 1}}$, which provides the bound we need for $|J_1(\lambda)| $ to give the correct exponent in part 2 of Theorem 1.1. On the other
hand if $\epsilon \geq 1$, then $(1.4)$ still holds for $\epsilon = 1$, in which case the second part of Lemma 3.1 gives a bound of $C'|\lambda|^{-{1 \over 2}}\ln|\lambda|$, which is 
better than $C''|\lambda|^{-\delta}$ for all $\delta < {1 \over 2}$, which is what we need for part 2 of Theorem 1.1 when $\epsilon \geq  1$. Thus we see that in all cases that 
$|J_1(\lambda)| $ is bounded by the appropriate quantity in part 2 of Theorem 1.1. This completes the proof of Theorem 1.1.

\section{References.}

\noindent [AGuV] V. Arnold, S Gusein-Zade, A Varchenko, {\it Singularities of differentiable maps
Volume II}, Birkhauser, Basel, 1988. \parskip = 4pt\baselineskip = 3pt

\noindent [BakMVaW] J. Bak, D. McMichael, J. Vance, S. Wainger, {\it Fourier transforms of surface area measure on convex surfaces in $\R^3$}, Amer. J. Math.{\bf 111} (1989), no.4,
 633-668. 

\noindent [BaGuZhZo] S. Basu, S. Guo, R. Zhang, P. Zorin-Kranich, {\it A stationary set method for estimating oscillatory integrals}, to appear, J. Eur. Math. Soc. 

\noindent [BrHoI] L. Brandolini, S. Hoffmann, A. Iosevich, {\it Sharp rate of average decay of the Fourier transform of a bounded set},
Geom. Funct. Anal. {\bf 13} (2003), no. 4, 671-680. 

\noindent [BNW] J. Bruna, A. Nagel, and S. Wainger, {\it Convex hypersurfaces and Fourier transforms},
Ann. of Math. (2) {\bf 127} no. 2, (1988), 333-365. 

\noindent [CaCWr] A. Carbery, M. Christ, J. Wright, {\it Multidimensional van der Corput and sublevel set estimates}, J. Amer. Math. Soc. {\bf 12}
 (1999), no. 4, 981-1015. 

\noindent [Gre] J. Green, {\it Lower bounds on $L^p$ quasi-norms and the uniform sublevel set problem}, Mathematika {\bf 67} (2021), no. 2, 296-323.

\noindent [G1] M. Greenblatt, {\it Oscillatory integral decay, sublevel set growth, and the Newton polyhedron}, Math. Annalen {\bf 346} (2010), no. 4, 857-895.

\noindent [G2] M. Greenblatt, {\it Maximal averages over hypersurfaces and the Newton polyhedron}, J. Funct. Anal. {\bf 262} (2012), no. 5, 2314-2348. 

\noindent [G3] M. Greenblatt, {\it A method for bounding oscillatory integrals in terms of non-oscillatory integrals}, submitted. 

\noindent [G4] M. Greenblatt, {\it Fourier transforms of irregular mixed homogeneous hypersurface measures},  Math. Nachr. {\bf 291} (2018), no.7, 1075-1087.

\noindent [G5] M. Greenblatt, {\it Fourier transforms of indicator functions, lattice point discrepancy, and the stability of integrals}, 
Math. Ann. {\bf 380} (2021), no.3-4, 1959-1990.

\noindent [Gr] P. Gressman, {\it Scalar oscillatory integrals in smooth spaces of homogeneous type}, 
Rev. Mat. Iberoam. {\bf 31} (2015), no. 1, 215–244. 

\noindent [H1] H. Hironaka, {\it Resolution of singularities of an algebraic variety over a field of characteristic zero I}, 
 Ann. of Math. (2) {\bf 79} (1964), 109-203.

\noindent [H2] H. Hironaka, {\it Resolution of singularities of an algebraic variety over a field of characteristic zero II},  
Ann. of Math. (2) {\bf 79} (1964), 205-326. 

\noindent [IkKeMu] I. Ikromov, M. Kempe, and D. M\"uller, {\it Estimates for maximal functions associated
to hypersurfaces in $R^3$ and related problems of harmonic analysis}, Acta Math. {\bf 204} (2010), no. 2,
151--271.

\noindent [IkMu] I. Ikromov, D. M\"uller, {\it  Uniform estimates for the Fourier transform of surface-carried measures in
$\R^3$ and an application to Fourier restriction}, J. Fourier Anal. Appl, {\bf 17} (2011), no. 6, 1292-1332.

\noindent [NaSeW] A. Nagel, A. Seeger, and S. Wainger, {\it Averages over convex hypersurfaces}, Amer. J. Math. {\bf 115} (1993), no. 4, 903-927.

\noindent [PhStS]  D. H. Phong, E. M. Stein, J. Sturm,  {\it On the growth and stability of real-analytic functions}, Amer. J. Math. {\bf 121} (1999), 
519-554.

\noindent [S] E. Stein, {\it Harmonic analysis; real-variable methods, orthogonality, and oscillatory \hfill\break
integrals}, Princeton Mathematics Series {\bf 43}, Princeton University Press, Princeton, NJ, 1993.

\noindent [V] A. N. Varchenko, {\it Newton polyhedra and estimates of oscillatory integrals}, Functional 
Anal. Appl. {\bf 18} (1976), no. 3, 175-196.

\vskip 0.5 in

\noindent Department of Mathematics, Statistics, and Computer Science \hfill \break
\noindent University of Illinois at Chicago \hfill \break
\noindent 322 Science and Engineering Offices \hfill \break
\noindent 851 S. Morgan Street \hfill \break
\noindent Chicago, IL 60607-7045 \hfill \break
\noindent greenbla@uic.edu

\end{document}